\documentclass[11pt,a4paper,twoside]{article}
\usepackage{latexsym,amsmath,amsthm,amssymb,amsfonts}
\usepackage{mathrsfs,multicol}
\usepackage{amscd}
\usepackage{cite}
\usepackage[ansinew]{inputenc}

\numberwithin{equation}{section}
\newcommand{\zdhy}{\allowdisplaybreak}

\catcode`@=11

\newskip\plaincentering \plaincentering=0pt plus 1000pt minus 1000pt
\def\@plainlign{\tabskip=0pt\everycr={}}
\def\eqalignno#1{\displ@y \tabskip\plaincentering
  \halign to\displaywidth{\hfil$\@lign\displaystyle{##}$\tabskip\z@skip
    &$\@lign\displaystyle{{}##}$\hfil\tabskip\plaincentering
    &\llap{$\@lign##$}\tabskip\z@skip\crcr
    #1\crcr}}
\def\leqalignno#1{\displ@y \tabskip\plaincentering
  \halign to\displaywidth{\hfil$\@lign\displaystyle{##}$\tabskip\z@skip
    &$\@lign\displaystyle{{}##}$\hfil\tabskip\plaincentering
    &\kern-\displaywidth\rlap{$\@lign##$}\tabskip\displaywidth\crcr
    #1\crcr}}
\def\plainLet@{\relax\iffalse{\fi\let\\=\cr\iffalse}\fi}
\def\plainvspace@{\def\vspace##1{\noalign{\vskip##1}}}

\def\intic@{\mathchoice{\hskip5\p@}{\hskip4\p@}{\hskip4\p@}{\hskip4\p@}}
\def\negintic@
 {\mathchoice{\hskip-5\p@}{\hskip-4\p@}{\hskip-4\p@}{\hskip-4\p@}}
\def\intkern@{\mathchoice{\!\!\!}{\!\!}{\!\!}{\!\!}}
\def\intdots@{\mathchoice{\cdots}{{\cdotp}\mkern1.5mu
    {\cdotp}\mkern1.5mu{\cdotp}}{{\cdotp}\mkern1mu{\cdotp}\mkern1mu
      {\cdotp}}{{\cdotp}\mkern1mu{\cdotp}\mkern1mu{\cdotp}}}
\newcount\intno@
\def\iint{\intno@=\tw@\futurelet\next\ints@}
\def\iiint{\intno@=\thr@@\futurelet\next\ints@}
\def\iiiint{\intno@=4 \futurelet\next\ints@}
\def\idotsint{\intno@=\z@\futurelet\next\ints@}
\def\ints@{\findlimits@\ints@@}
\newif\iflimtoken@
\newif\iflimits@
\def\findlimits@{\limtoken@false\limits@false\ifx\next\limits
 \limtoken@true\limits@true\else\ifx\next\nolimits\limtoken@true\limits@false
    \fi\fi}
\def\multintlimits@{\intop\ifnum\intno@=\z@\intdots@
  \else\intkern@\fi
    \ifnum\intno@>\tw@\intop\intkern@\fi
     \ifnum\intno@>\thr@@\intop\intkern@\fi\intop}
\def\multint@{\int\ifnum\intno@=\z@\intdots@\else\intkern@\fi
   \ifnum\intno@>\tw@\int\intkern@\fi
    \ifnum\intno@>\thr@@\int\intkern@\fi\int}
\def\ints@@{\iflimtoken@\def\ints@@@{\iflimits@
   \negintic@\mathop{\intic@\multintlimits@}\limits\else
    \multint@\nolimits\fi\eat@}\else
     \def\ints@@@{\multint@\nolimits}\fi\ints@@@}
\def\Sb{_\bgroup\vspace@
        \baselineskip=\fontdimen10 \scriptfont\tw@
        \advance\baselineskip by \fontdimen12 \scriptfont\tw@
        \lineskip=\thr@@\fontdimen8 \scriptfont\thr@@
        \lineskiplimit=\thr@@\fontdimen8 \scriptfont\thr@@
        \Let@\vbox\bgroup\halign\bgroup \hfil$\scriptstyle
            {##}$\hfil\cr}
\def\endSb{\crcr\egroup\egroup\egroup}
\def\Sp{^\bgroup\vspace@
        \baselineskip=\fontdimen10 \scriptfont\tw@
        \advance\baselineskip by \fontdimen12 \scriptfont\tw@
        \lineskip=\thr@@\fontdimen8 \scriptfont\thr@@
        \lineskiplimit=\thr@@\fontdimen8 \scriptfont\thr@@
        \Let@\vbox\bgroup\halign\bgroup \hfil$\scriptstyle
            {##}$\hfil\cr}
\def\endSp{\crcr\egroup\egroup\egroup}
\def\Let@{\relax\iffalse{\fi\let\\=\cr\iffalse}\fi}
\def\vspace@{\def\vspace##1{\noalign{\vskip##1 }}}
\def\aligned{\,\vcenter\bgroup\plainvspace@\plainLet@\openup\jot\m@th\ialign
  \bgroup \strut\hfil$\displaystyle{##}$&$\displaystyle{{}##}$\hfil\crcr}
\def\endaligned{\crcr\egroup\egroup}
\def\matrix{\,\vcenter\bgroup\plainLet@\plainvspace@
    \normalbaselines
  \m@th\ialign\bgroup\hfil$##$\hfil&&\quad\hfil$##$\hfil\crcr
    \mathstrut\crcr\noalign{\kern-\baselineskip}}
\def\endmatrix{\crcr\mathstrut\crcr\noalign{\kern-\baselineskip}\egroup
                \egroup\,}
\newtoks\hashtoks@
\hashtoks@={#}
\def\format{\crcr\egroup\iffalse{\fi\ifnum`}=0 \fi\format@}
\def\format@#1\\{\def\preamble@{#1}%
  \def\c{\hfil$\the\hashtoks@$\hfil}%
  \def\r{\hfil$\the\hashtoks@$}%
  \def\l{$\the\hashtoks@$\hfil}%
  \setbox\z@=\hbox{\xdef\Preamble@{\preamble@}}\ifnum`{=0 \fi\iffalse}\fi
   \ialign\bgroup\span\Preamble@\crcr}

\def\cases{\left\{\,\vcenter\bgroup\plainvspace@
     \normalbaselines\openup\jot\m@th
      \plainLet@\ialign\bgroup$\displaystyle{##}$\hfil&
      \quad$\displaystyle{{}##}$\hfil\crcr
      \mathstrut\crcr\noalign{\kern-\baselineskip}}

\newif\iftagsleft@
\tagsleft@true
\def\TagsOnRight{\global\tagsleft@false}
\def\tag#1$${\iftagsleft@\leqno\else\eqno\fi
 \hbox{\def\pagebreak{\global\postdisplaypenalty-\@M}%
 \def\nopagebreak{\global\postdisplaypenalty\@M}\rm(#1\unskip)}%
  $$\postdisplaypenalty\z@\ignorespaces}
\interdisplaylinepenalty=\@M
\def\allowdisplaybreak{\noalign{\allowbreak}}
\def\plainallowdisplaybreak@{\def\allowdisplaybreak{\noalign{\allowbreak}}}
\def\plaindisplaybreak@{\def\displaybreak{\noalign{\break}}}
\def\align#1\endalign{\def\tag{&}\plainvspace@\plainallowdisplaybreak@
\plaindisplaybreak@
  \iftagsleft@\plainlalign@#1\endalign\else
   \plainralign@#1\endalign\fi}
\def\plainralign@#1\endalign{\displ@y\plainLet@\tabskip\plaincentering
\halign to\displaywidth
     {\hfil$\displaystyle{##}$\tabskip=\z@&$\displaystyle{{}##}$\hfil
       \tabskip=\plaincentering&\llap{\hbox{\rm(##\unskip)}}\tabskip\z@\crcr
             #1\crcr}}
\def\plainlalign@
 #1\endalign{\displ@y\plainLet@\tabskip\plaincentering\halign to \displaywidth
   {\hfil$\displaystyle{##}$\tabskip=\z@&$\displaystyle{{}##}$\hfil
   \tabskip=\plaincentering&\kern-\displaywidth
        \rlap{\hbox{\rm(##\unskip)}}\tabskip=\displaywidth\crcr
               #1\crcr}}

\def\re@#1{\par\hangindent\parindent\indent\llap{#1\enspace}\ignorespaces}
\def\qfootnote#1{\edef\@sf{\spacefactor\the\spacefactor}{}#1\@sf
      \insert\footins{\let\egroup=}\footnotesize 
      \interlinepenalty100 \let\par=\endgraf
        \leftskip=0pt \rightskip=0pt
        \splittopskip=10pt plus 1pt minus 1pt \floatingpenalty=20000
   \smallskip\re@{#1}\bgroup\strut\aftergroup{\strut\egroup}\let\next}
\catcode`\@=\active
\TagsOnRight
\begin{document}
\title{\bf Gradient estimates for the porous medium
equations on Riemannian manifolds\footnote{The research of the first
author is supported by NSFC grant(No. 11001076) and NSF of Henan
Provincial Education department(No. 2010A110008). The research of
the third author is supported by NSFC grant(No. 10971110).}}
\author{Guangyue Huang,\ Zhijie Huang,\ Haizhong Li
}
\date{}
\maketitle
\begin{quotation}
\noindent{\bf Abstract.}~ In this paper we study gradient estimates
for the positive solutions of the porous medium equation:
$$u_t=\Delta u^m$$
where $m>1$, which is a nonlinear version of the heat equation. We
derive local gradient estimates of the Li-Yau type for positive
solutions of porous medium equations on Riemannian manifolds with
Ricci curvature bounded from below. As applications, several
parabolic Harnack inequalities are obtained. In particular, our
results improve the ones of Lu, Ni, V\'{a}zquez and Villani in
\cite{pengni09}. Moreover, our results recover the
ones of Davies in \cite{davies89}, Hamilton in \cite{hamilton93} and Li and Xu in \cite{lixu11}.
\end{quotation}

\section{Introduction}

Let $(M^n,g)$ be an $n$-dimensional complete Riemannian manifold. Li
and Yau \cite{{liyau86}} studied positive solutions of the heat
equation
\begin{equation}\label{maysec1}
u_t=\Delta u
\end{equation} and obtained the following gradient estimates:

\noindent{\bf Theorem A(Li-Yau\cite{liyau86}).} {\it Let $(M^n,g)$
be a complete Riemannian manifold with ${\rm Ric}(B_p(2R))\\
\geq-K$, $K\geq0$. Suppose that $u$ is a positive solution of
\eqref{maysec1} on $B_p(2R)\times [0,T]$. Then on $B_p(R)$, we have
\begin{equation}\label{maysec2}
\frac{|\nabla u|^2}{u^2} - \alpha\frac{u_t}{u}
\leq\frac{C\alpha^2}{R^2}\left(\frac{\alpha^2}{\alpha-1}+\sqrt{K} R
\right)+\frac{n\alpha^2K}{2(\alpha-1)}+\frac{n\alpha^2}{2t},
\end{equation} where $\alpha>1$ is a constant and the constant $C$
depends only on $n$. Moreover, when $R\rightarrow \infty$,
\eqref{maysec2} yields the following estimate on complete noncompact
Riemannian manifold $(M^n,g)$:
\begin{equation}\label{maysec3}
\frac{|\nabla u|^2}{u^2}-\alpha\frac{u_t}{u} \leq
\frac{n\alpha^2K}{2(\alpha-1)}+\frac{n\alpha^2}{2t}.\end{equation} }

In \cite{davies89}, Davies improved the estimate \eqref{maysec3} to

\noindent{\bf Theorem B(Davies\cite{davies89}).} {\it Let $(M^n,g)$
be a complete noncompact Riemannian manifold with ${\rm Ric}\geq-K$,
$K\geq0$. Suppose that $u$ is a positive solution to
\eqref{maysec1}. Then
\begin{equation}\label{maysec4}
\frac{|\nabla u|^2}{u^2}-\alpha\frac{u_t}{u} \leq
\frac{n\alpha^2K}{4(\alpha-1)}+\frac{n\alpha^2}{2t}.\end{equation}

}

In \cite{hamilton93}, Hamilton proved the following estimate:

\noindent{\bf Theorem C(Hamilton\cite{hamilton93}).} {\it Let
$(M^n,g)$ be a complete noncompact Riemannian manifold with ${\rm
Ric}\geq-K$, $K\geq0$. Suppose that $u$ is a positive solution to
\eqref{maysec1}. Then
\begin{equation}\label{hamilton}
\frac{|\nabla u|^2}{u^2}-e^{2Kt}\frac{u_t}{u} \leq
e^{4Kt}\frac{n}{2t}.\end{equation}

}

Recently, Li and Xu \cite{lixu11} obtained new Li-Yau type gradient
estimates for positive solutions of the heat equation
\eqref{maysec1} on Riemannian manifolds. For the related research
and improvement in this direction, see
\cite{Yau94,Yau95,Bakry99,hamilton93,Ma06,Chen09,Lixiangdong,Yangyuyan07,huangma10}
and the references therein.

The porous medium equation
\begin{equation}\label{sec1}
u_t=\Delta u^m,
\end{equation} where $m>1$ is a nonlinear version of the heat
equation \eqref{maysec1}. For various values of $m>1$, it has arisen
in different applications to model diffusive phenomena. The readers
who are interested in the applications of \eqref{sec1} see
\cite{pengni09,vazquez07,Aronson79} and the references therein. In
\cite{pengni09}, Lu, Ni, V\'{a}zquez and Villani studied gradient
estimates of \eqref{sec1} with $m>1$ and proved the following
results (see Theorem 3.3 in \cite{pengni09}):

\noindent{\bf Theorem D(P. Lu, L. Ni, J. V\'{a}zquez,
C.Villani\cite{pengni09}).} {\it Let $(M^n,g)$ be a complete
Riemannian manifold with ${\rm Ric}(B_p(2R))\geq-K$, $K\geq0$.
Suppose that $u$ is a positive solution to \eqref{sec1}. Let
$v=\frac{m}{m-1}u^{m-1}$ and $M=\max_{B_p(2R)\times [0,T]}v$. Then
for any $\alpha>1$, on the ball $B_p(R)$, we have
\begin{equation}\label{maysec5}\aligned
 \frac{|\nabla v|^2}{v} -
\alpha\frac{v_t}{v}
\leq&\frac{CMa\alpha^2}{R^2}\left(\frac{\alpha^2}{\alpha-1}am^2+(m-1)(1+\sqrt{K}R)
\right)\\
&+\frac{\alpha^2}{\alpha-1}a(m-1)MK+\frac{a\alpha^2}{t},
\endaligned\end{equation}
where $a=\frac{n(m-1)}{n(m-1)+2}$ and the constant $C$ depends only
on $n$. Moreover, when $R\rightarrow \infty$, \eqref{maysec5} yields
the following estimate on complete noncompact Riemannian manifold
$(M^n,g)$:\begin{equation}\label{maysec5-1}\aligned
 \frac{|\nabla v|^2}{v} -
\alpha\frac{v_t}{v}\leq&\frac{\alpha^2}{\alpha-1}a(m-1)MK+
\frac{a\alpha^2}{t}.
\endaligned\end{equation}}

In this paper, we further study gradient estimates of the porous
medium equation \eqref{sec1}.  We derive Davies's type estimate and
Hamilton's
 to \eqref{sec1}. Besides, we obtain estimates of Li-Xu
 type for \eqref{sec1}. In particular, our results improve the ones
 of Lu, Ni, V\'{a}zquez and Villani in \cite{pengni09}. Now, we state our
results as follows:

\noindent{\bf Theorem 1.1.} {\it Let $(M^n,g)$ be a complete
Riemannian manifold with ${\rm Ric}(B_p(2R))\geq-K$, $K\geq0$.
Suppose that $u$ is a positive solution to \eqref{sec1}. Let
$v=\frac{m}{m-1}u^{m-1}$ and $M=\max_{B_p(2R)\times [0,T]}v$. Then
for any $\alpha>1$, on the ball $B_p(R)$, we have
\begin{equation}\label{thmhuang1}\aligned
 \frac{|\nabla v|^2}{v} -
\alpha\frac{v_t}{v}\leq&\left\{\frac{ a \alpha^2m M^{\frac{1}{2}}
}{(\alpha-1)^{\frac{1}{2}}}\frac{C}{R}+a^{\frac{1}{2}}\alpha \left[
\frac{1}{t} +\frac{(m-1)MK}{2(\alpha -1)} \right. \right. \\
&\left.\left. +(m-1)M\frac{C}{R^2}\left(1+ \sqrt{K} \coth (\sqrt{K}
R)\right) \right]^{\frac{1}{2}}\right\}^2,
\endaligned\end{equation}
where $a=\frac{n(m-1)}{n(m-1)+2}$ and the constant $C$ depends only
on $n$.
}

Letting $R\rightarrow \infty$, we obtain the gradient estimates on
complete noncompact Riemannian manifolds, which improves
\eqref{maysec5-1} of Theorem D in \cite{pengni09}.

\noindent{\bf Corollary 1.1.} {\it Let $(M^n,g)$ be a complete
noncompact Riemannian manifold with ${\rm Ric}\geq-K$, $K\geq0$.
Suppose that $u$ is a positive solution to \eqref{sec1}. Let
$v=\frac{m}{m-1}u^{m-1}$ and $M=\sup_{M^n\times [0,T]}v$. Then for
any $\alpha>1$, we have
\begin{equation}\label{thmhuang2}\aligned
 \frac{|\nabla v|^2}{v} -
\alpha\frac{v_t}{v}\leq&\frac{\alpha^2}{2(\alpha-1)}a(m-1)MK+
\frac{a\alpha^2}{t}.
\endaligned\end{equation} }

Applying the inequality \eqref{thmhuang2}, we derive the following
Harnack inequality:

\noindent{\bf Corollary 1.2.} {\it Let $(M^n,g)$ be a complete
noncompact Riemannian manifold with ${\rm Ric} \geq-K$, $K\geq0$.
Suppose that $u$ is a positive solution to \eqref{sec1}. Let
$v=\frac{m}{m-1}u^{m-1}$ and $M=\sup_{M^n\times [0,T]}v,\,
\tilde{M}=\inf_{M^n\times [0,T]}v$. Then for any $x_1,x_2\in M^n$,
$0<t_1<t_2<T$, $\alpha>1$, the following inequality holds:
\begin{equation}\label{thmhuang3}\aligned
v(x_1,t_1)\leq&v(x_2,t_2)\left(\frac{t_2}{t_1} \right)^{a\alpha}\exp
\left\{\frac{\alpha\,{dist}^2(x_2,x_1)}{4\tilde{M}(t_2-t_1)}\right.\\
&\left. +\frac{\alpha}{2(\alpha-1)}a(m-1)MK(t_2-t_1) \right\},
\endaligned\end{equation} where
${dist}(x_2,x_1)$ is the distance between $x_1$ and $x_2$.
 }

\noindent{\bf Remark 1.1.}  We rewrite the inequality
\eqref{maysec5-1} as
\begin{equation}\label{remark11}\aligned
|\nabla v|^2  - \alpha v_t\leq&\frac{\alpha^2}{\alpha-1}a(m-1)MKv+
\frac{a\alpha^2v}{t}.
\endaligned\end{equation}
Since $(m-1)v=mu^{m-1}$, we have $(m-1)v\rightarrow 1$ as
$m\rightarrow 1$. Hence, $(m-1)M\rightarrow 1$,
$$\aligned
|\nabla v|^2&\rightarrow \frac{|\nabla
u|^2}{u^2},\\
 v_t&\rightarrow \frac{u_t}{u},\\
 av &\rightarrow \frac{n}{2}, \endaligned$$ as $m\rightarrow 1$.
As a result, \eqref{remark11} becomes the inequality \eqref{maysec3}
in Theorem A of Li-Yau. Therefore, for complete noncompact
Riemannian manifold $(M^n,g)$, the estimate \eqref{maysec5-1} in
Theorem D of Lu, Ni, V\'{a}zquez and Villani reduces to the estimate
\eqref{maysec3} in Theorem A of Li-Yau when $m\rightarrow 1$.
Similarly, it is easy to see that our estimates \eqref{thmhuang2}
reduces to the estimate \eqref{maysec4} in Theorem B of Davies when
$m\rightarrow 1$.

\noindent{\bf Theorem 1.2.} {\it Let $(M^n,g)$ be a complete
Riemannian manifold with ${\rm Ric}(B_p(2R))\geq-K$, $K\geq0$.
Suppose that $u$ is a positive solution to \eqref{sec1}. Let
$v=\frac{m}{m-1}u^{m-1}$ and $M=\max_{B_p(2R)\times [0,T]}v$. Then
on the ball $B_p(R)$, we have
\begin{alignat}{1}\label{thhamilton1}
 \frac{|\nabla v|^2}{v}-\alpha(t)\frac{v_t}{v}
 \leq&\left(\frac{m^2Ma^2\alpha^4(t)}{2(\alpha(t)-1)}
+3(m-1)Ma\alpha^2(t)\right)\frac{C}{R^2}\\
&\quad+(m-1)Ma\alpha^2(t)\sqrt{K}\coth(\sqrt{K}R)\frac{C}{R}+\frac{a\alpha^2(t)}{t},
\end{alignat}
where $a=\frac{n(m-1)}{n(m-1)+2}$, $\alpha(t)=e^{2(m-1)MKt}$ and the
constant $C$ depends only on $n$. }

Letting $R\rightarrow \infty$, we obtain the following gradient
estimates on complete noncompact Riemannian manifolds.

\noindent{\bf Corollary 1.3.} {\it Let $(M^n,g)$ be a complete
noncompact Riemannian manifold with ${\rm Ric}\geq-K$, $K\geq0$.
Suppose that $u$ is a positive solution to \eqref{sec1}. Let
$v=\frac{m}{m-1}u^{m-1}$ and $M=\sup_{M^n\times [0,T]}v$. Then we
have
\begin{equation}\label{thhamilton2}\aligned
 \frac{|\nabla v|^2}{v} -
\alpha(t)\frac{v_t}{v}\leq&\frac{a\alpha^2(t)}{t},
\endaligned\end{equation} where $\alpha(t)=e^{2(m-1)MKt}$.}

Applying the inequality \eqref{thhamilton2}, we derive the following
Harnack inequality:

\noindent{\bf Corollary 1.4.} {\it Let $(M^n,g)$ be a complete
noncompact Riemannian manifold with ${\rm Ric} \geq-K$, $K\geq0$.
Suppose that $u$ is a positive solution to \eqref{sec1}. Let
$v=\frac{m}{m-1}u^{m-1}$ and $M=\sup_{M^n\times [0,T]}v,\,
\tilde{M}=\inf_{M^n\times [0,T]}v$. Then for any $x_1,x_2\in M^n$,
$0<t_1<t_2<T$, $\alpha>1$, the following inequality holds:
\begin{alignat}{1}\label{thhamilton3}
v(x_1,t_1)\leq&v(x_2,t_2)\exp\left\{\frac{e^{2(m-1)MKt_2}-e^{2(m-1)MKt_1}}{2(m-1)MK}\left(
\frac{{dist}^2(x_2,x_1)}{4\tilde{M}(t_2-t_1)^2}
+\frac{a}{t_1}\right)\right\},
\end{alignat} where
${dist}(x_2,x_1)$ is the distance between $x_1$ and $x_2$.
 }

\noindent{\bf Remark 1.2.}  Notice that \eqref{thhamilton2} can be
written as
\begin{equation}\label{remarkhamilton1}\aligned
|\nabla v|^2-\alpha(t) v_t\leq&\frac{av}{t}\alpha^2(t).
\endaligned\end{equation}
Since $(m-1)v=mu^{m-1}$, we have $(m-1)v\rightarrow 1$ as
$m\rightarrow 1$. Hence, $(m-1)M\rightarrow 1$,
$av\rightarrow\frac{n}{2}$ and $\alpha(t)\rightarrow e^{2Kt}$. Hence
letting $m\rightarrow 1$ in \eqref{remarkhamilton1} yields the
inequality \eqref{hamilton}. Therefore, our Corollary 1.3 extends
 Theorem C of Hamilton.

\noindent{\bf Theorem 1.3.} {\it Let $(M^n,g)$ be a complete
Riemannian manifold with ${\rm Ric}(B_p(2R))\geq-K$, $K\geq0$.
Suppose that $u$ is a positive solution to \eqref{sec1}. Let
$v=\frac{m}{m-1}u^{m-1}$ and $M=\max_{B_p(2R)\times [0,T]}v$. Then
on the ball $B_p(R)$, we have
\begin{equation}\label{maysec6}\aligned
 \frac{|\nabla v|^2}{v} -
\alpha(t)\frac{v_t}{v}-\varphi(t) \leq&\left\{a(m-1)
\left(\frac{C}{R^2}+\frac{C\sqrt{K}\coth(\sqrt{K}R)}{R}\right)\right.\\
&\left.+a^2m^2\frac{C}{R^2\tanh((m-1)MKt)}\right\}M,
\endaligned\end{equation}
where $a=\frac{n(m-1)}{n(m-1)+2}$ and the constant $C$ depends only
on $n$. $\alpha(t)$ and $\varphi(t)$ are given by
\begin{alignat}{1}\label{maysec7}
\varphi(t)
=&a(m-1)MK\{\coth((m-1)MKt)+1\},\\
\alpha(t)=&1+\frac{\cosh((m-1)MKt)\sinh((m-1)MKt)-(m-1)MKt}{\sinh^2((m-1)MKt)}.
\end{alignat}
}

Letting $R\rightarrow \infty$, we obtain the gradient estimates on
complete noncompact Riemannian manifolds.

\noindent{\bf Corollary 1.5.} {\it Let $(M^n,g)$ be a complete
noncompact Riemannian manifold with ${\rm Ric}\geq-K$, $K\geq0$.
Suppose that $u$ is a positive solution to \eqref{sec1}. Let
$v=\frac{m}{m-1}u^{m-1}$ and $M=\sup_{M^n\times [0,T]}v$. Then we
have
\begin{equation}\label{maysec8}
 \frac{|\nabla v|^2}{v} -
\alpha(t)\frac{v_t}{v}-\varphi(t) \leq0,
\end{equation}
where $\alpha(t)$ and $\varphi(t)$ are given by \eqref{maysec7}. }

Applying the inequality \eqref{maysec8}, we derive the following
Harnack inequality:

\noindent{\bf Corollary 1.6.} {\it Let $(M^n,g)$ be a complete
noncompact Riemannian manifold with ${\rm Ric} \geq-K$. Suppose that
$u$ is a positive solution to \eqref{sec1}. Let
$v=\frac{m}{m-1}u^{m-1}$ and $M=\sup_{M^n\times [0,T]}v,\,
\tilde{M}=\inf_{M^n\times [0,T]}v$. Then for any $x_1,x_2\in M^n$,
$0<t_1<t_2<T$, the following inequality holds:
\begin{equation}\label{maysec9}
v(x_1,t_1)\leq v(x_2,t_2)A_1(t_1,t_2)
\exp\left\{\frac{{dist}^2(x_2,x_1)}{4\tilde{M}(t_2-t_1)}
(1+A_2(t_1,t_2))\right\},
\end{equation}
where ${dist}(x_2,x_1)$ is the distance between $x_1$ and $x_2$.
Moreover,
$$A_1(t_1,t_2)=\left(
\frac{\exp(2(m-1)MKt_2)-2(m-1)MKt_2-1}{\exp(2(m-1)MKt_1)
-2(m-1)MKt_1-1}\right)^{\frac{a}{2}},$$
$$A_2(t_1,t_2)=\frac{t_2\coth((m-1)MKt_2)-t_1\coth((m-1)MKt_1)}{t_2-t_1}.$$
 }

A linear version of Theorem 1.3 is the following:

\noindent{\bf Theorem 1.4.} {\it Let $(M^n,g)$ be a complete
Riemannian manifold with ${\rm Ric}(B_p(2R))\geq-K$, $K\geq0$.
Suppose that $u$ is a positive solution to \eqref{sec1}. Let
$v=\frac{m}{m-1}u^{m-1}$ and $M=\max_{B_p(2R)\times [0,T]}v$. Then
on the ball $B_p(R)$, we have
\begin{alignat}{1}\label{maysec10}
 \frac{|\nabla v|^2}{v} -
\alpha(t)\frac{v_t}{v}-\varphi(t) \leq&\left\{a(m-1)\alpha^2(t)
\left(
\frac{C}{R^2}+\frac{C\sqrt{K}\coth(\sqrt{K}R)}{R}\right)\right.\\
& \left.+\frac{a^2m^2\alpha^4(t)}{\beta(t)}\frac{C}{R^2}\right\}M,
\end{alignat}
where $a=\frac{n(m-1)}{n(m-1)+2}$ and the constant $C$ depends only
on $n$. $\alpha(t)$ and $\varphi(t)$ are given by
\begin{equation}\label{maysec11}\aligned \varphi(t)
=&\frac{a}{t}+a(m-1)MK+\frac{a}{3}((m-1)MK)^2t,\\
\alpha(t)=&1+\frac{2}{3}(m-1)MKt.
\endaligned\end{equation}
}

Letting $R\rightarrow \infty$, we obtain the gradient estimates on
complete noncompact Riemannian manifolds.

\noindent{\bf Corollary 1.7.} {\it Let $(M^n,g)$ be a complete
noncompact Riemannian manifold with ${\rm Ric}\geq-K$, $K\geq0$..
Suppose that $u$ is a positive solution to \eqref{sec1}. Let
$v=\frac{m}{m-1}u^{m-1}$ and $M=\sup_{M^n\times [0,T]}v$. Then we
have
\begin{equation}\label{maysec12}
 \frac{|\nabla v|^2}{v} -
\alpha(t)\frac{v_t}{v}-\varphi(t) \leq0,
\end{equation}
where $\alpha(t)$ and $\varphi(t)$ are given by \eqref{maysec11}. }

Applying the inequality \eqref{maysec12}, we derive the following
Harnack inequality:

\noindent{\bf Corollary 1.8.} {\it Let $(M^n,g)$ be a complete
noncompact Riemannian manifold with ${\rm Ric} \geq-K$, $K\geq0$.
Suppose that $u$ is a positive solution to \eqref{sec1}. Let
$v=\frac{m}{m-1}u^{m-1}$ and $M=\sup_{M^n\times [0,T]}v,\,
\tilde{M}=\inf_{M^n\times [0,T]}v$. Then for any $x_1,x_2\in M^n$,
$0<t_1<t_2<T$, the following inequality holds:
\begin{equation}\label{maysec13}\aligned
v(x_1,t_1)\leq& v(x_2,t_2)\left(\frac{t_2}{t_1}\right)^a
\left(\frac{1+\frac{2}{3}(m-1)MKt_2}{1+\frac{2}{3}(m-1)MKt_1}
\right)^{\frac{-a}{4}}\\
&\exp\left\{\frac{{dist}^2(x_2,x_1)}{4\tilde{M}(t_2-t_1)}\left(1+
\frac{1}{3}(m-1)MK(t_2+t_1)\right)+\frac{a}{2}(m-1)MK(t_2-t_1)\right\},
\endaligned\end{equation} where
${dist}(x_2,x_1)$ is the distance between $x_1$ and $x_2$.
 }

\noindent{\bf Remark 1.3.} When $m\rightarrow 1$, our Theorem 1.3
reduces to Theorem 1.1 of Li and Xu in \cite{lixu11}. Similarly, our
Theorem 1.4 reduces to Theorem 1.2 of Li and Xu in \cite{lixu11}.
Note that \eqref{maysec6} can be written as
\begin{equation}\label{remark1}\aligned
|\nabla v|^2- \alpha(t)v_t-\varphi(t)v \leq&\left\{a(m-1)
\left(\frac{C}{R^2}+\frac{C\sqrt{K}\coth(\sqrt{K}R)}{R}\right)\right.\\
&\left.+a^2m^2\frac{C}{R^2\tanh((m-1)MKt)}\right\}Mv.
\endaligned\end{equation}
Since $(m-1)v=mu^{m-1}$, we have $(m-1)v\rightarrow 1$ as
$m\rightarrow 1$. Hence, $(m-1)M\rightarrow 1$,
\begin{alignat*}{1}
&\alpha(t)\rightarrow 1+\frac{\cosh(Kt)\sinh(Kt)-Kt}{\sinh^2(Kt)},\\
\zdhy &\varphi(t)v\rightarrow\frac{nK}{2}\{\coth(Kt)+1\},\\ \zdhy
&|\nabla v|^2\rightarrow \frac{|\nabla u|^2}{u^2},\\  \zdhy &
v_t\rightarrow \frac{u_t}{u},\end{alignat*} and
$$\aligned
&\left\{a(m-1)
\left(\frac{C}{R^2}+\frac{C\sqrt{K}\coth(\sqrt{K}R)}{R}\right)\right.\\
&\left.+a^2m^2\frac{C}{R^2\tanh((m-1)MKt)}\right\}Mv
\rightarrow\frac{nC}{R^2}+\frac{nC\sqrt{K}}{R}\coth(\sqrt{K}R)+
\frac{n^2C}{R^2\tanh(Kt)}
\endaligned$$
as $m\rightarrow 1$. As a result, \eqref{remark1} becomes
\begin{alignat}{1}\label{remark2}
\frac{|\nabla
u|^2}{u^2}-\tilde{\alpha}(t)\frac{u_t}{u}-\tilde{\varphi}(t)\leq
\frac{nC}{R^2}+\frac{nC\sqrt{K}}{R}\coth(\sqrt{K}R)+
\frac{n^2C}{R^2\tanh(Kt)}
\end{alignat}
 by letting $m\rightarrow 1$, where $\tilde{\alpha}(t),\tilde{\varphi}(t)$ in
 \eqref{remark2} are given by
 $\tilde{\alpha}(t)=1+\frac{\cosh(Kt)\sinh(Kt)-Kt}{\sinh^2(Kt)}$,
 $\tilde{\varphi}(t)=\frac{nK}{2}\{\coth(Kt)+1\}$. Therefore, our Theorem
 1.3 becomes Theorem 1.1 of Li and Xu in
 \cite{lixu11} as long as letting $m\rightarrow 1$. Similarly, our Theorem
 1.4 becomes Theorem 1.2 of Li and Xu in
 \cite{lixu11} as long as letting $m\rightarrow 1$.

\noindent{\bf Remark 1.4.} When $t$ is small enough,
$\alpha(t),\varphi(t)$ defined by \eqref{maysec7} and
\eqref{maysec11} both satisfy $\alpha(t)\rightarrow 1$ and
$\varphi(t)\leq2a(m-1)MK+\frac{a}{t}$. Hence, by Corollary 1.5 and
1.7,  we have \begin{equation}\label{remark4}
 \frac{|\nabla v|^2}{v} -
\alpha(t)\frac{v_t}{v} \leq 2a(m-1)MK+\frac{a}{t}.
\end{equation} Clearly, for $t$ small enough, \eqref{remark4} is better than
\eqref{maysec5-1}. Thus Corollary 1.5 and 1.7 improve
\eqref{maysec5-1} in Theorem D of \cite{pengni09} in this sense.

\section{Proof of Theorem 1.1}

Let $v=\frac{m}{m-1}u^{m-1}$. From the equation \eqref{sec1}, one
gets $v_t=(m-1)v\Delta v+|\nabla v|^2$ which is equivalent to the
following form:
\begin{equation}\label{sec2}
\frac{v_t}{v}=(m-1)\Delta v+\frac{|\nabla v|^2}{v}.
\end{equation}

\noindent{\bf Lemma 2.1.} {\it As in \cite{pengni09}, we introduce
the differential operator
$$\mathcal{L}=\partial_t-(m-1)v\Delta.$$ Denote by $F=\frac{|\nabla
v|^2}{v}-\alpha\frac{v_t}{v}-\varphi$, where $\alpha=\alpha(t)$ and
$\varphi=\varphi(t)$ are functions depending on $t$, then we have
\begin{equation}\label{sec9}\aligned
\mathcal{L}(F)=&-2(m-1)v_{ij}^2-2(m-1)R_{ij}v_iv_j+2m \nabla v\nabla F\\
&-((m-1)\Delta v)^2+(1-\alpha)\left(\frac{v_t}{v}\right)^2
-\alpha'\frac{v_t}{v}-\varphi'.
\endaligned\end{equation}
}

\noindent{\it Proof of Lemma 2.1.} We need two formulas (see p5-p6
in \cite{pengni09})
$$\mathcal{L}\left(\frac{v_t}{v}\right)=(m-1)\Delta
v\frac{v_t}{v}+\frac{2}{v}\nabla v\nabla
v_t-\frac{v_t}{v}\frac{|\nabla
v|^2}{v}+2(m-1)v\nabla\left(\frac{v_t}{v}\right)\nabla(\log v),
$$
\begin{alignat*}{1} \mathcal{L}\left(\frac{|\nabla
v|^2}{v}\right)=&2(m-1)\Delta v\frac{|\nabla
v|^2}{v}+\frac{2}{v}\nabla v\nabla|\nabla
v|^2-2(m-1)v_{ij}^2\\
&-2(m-1)R_{ij}v_iv_j-\frac{|\nabla
v|^4}{v^2}+2(m-1)v\nabla\left(\frac{|\nabla
v|^2}{v}\right)\nabla(\log v).\end{alignat*} Hence we have
\begin{alignat}{1}\label{sec5}
\mathcal{L}(F)=&\mathcal{L}\left(\frac{|\nabla
v|^2}{v}\right)-\alpha\mathcal{L}\left(\frac{v_t}{v}\right)
-\alpha'\frac{v_t}{v}-\varphi'\\
=&2(m-1)\Delta v\frac{|\nabla v|^2}{v}+\frac{2}{v}\nabla
v\nabla|\nabla
v|^2-2(m-1)v_{ij}^2\\
&-2(m-1)R_{ij}v_iv_j-\frac{|\nabla v|^4}{v^2}
+2(m-1)v\nabla\left(\frac{|\nabla v|^2}{v}\right)\nabla(\log v)\\
&-\alpha(m-1)\Delta v\frac{v_t}{v}-\alpha\frac{2}{v}\nabla v\nabla
v_t+\alpha\frac{v_t}{v}\frac{|\nabla
v|^2}{v}-2\alpha(m-1)v\nabla\left(\frac{v_t}{v}\right)\nabla(\log v)\\
&-\alpha'\frac{v_t}{v}-\varphi'.
\end{alignat}
Since $$2(m-1)v\nabla(\frac{|\nabla v|^2}{v})\nabla(\log
v)-2\alpha(m-1)v\nabla(\frac{v_t}{v})\nabla(\log v)=2(m-1) \nabla
v\nabla F,
$$
$$\frac{2}{v}\nabla v\nabla|\nabla v|^2-\alpha\frac{2}{v}\nabla
v\nabla v_t=\frac{2}{v}\nabla v\nabla((F+\varphi)v)
=2(F+\varphi)\frac{|\nabla v|^2}{v}+2\nabla v\nabla F,
$$
we have
\begin{equation}\label{sec6}\aligned
2(m-1)&v\nabla\left(\frac{|\nabla v|^2}{v}\right)\nabla(\log
v)-2\alpha(m-1)v\nabla\left(\frac{v_t}{v}\right)\nabla(\log
v)+\frac{2}{v}\nabla v\nabla|\nabla v|^2-\alpha\frac{2}{v}\nabla
v\nabla v_t\\
=&2m \nabla v\nabla F+2(F+\varphi)\frac{|\nabla v|^2}{v} \\
=& 2m \nabla v\nabla F+2\left(\frac{|\nabla
v|^2}{v}-\alpha\frac{v_t}{v}\right)\frac{|\nabla v|^2}{v}.
\endaligned\end{equation}
On the other hand, it follows from \eqref{sec2} that
\begin{alignat}{1}\label{sec7}
2(m-1)&\Delta v\frac{|\nabla v|^2}{v}-\frac{|\nabla v|^4}{v^2}
-\alpha(m-1)\Delta v\frac{v_t}{v}+\alpha\frac{v_t}{v}\frac{|\nabla
v|^2}{v}\\
=&2\frac{|\nabla v|^2}{v}\left(\frac{v_t}{v}-\frac{|\nabla
v|^2}{v}\right)-\frac{|\nabla
v|^4}{v^2}-\alpha\frac{v_t}{v}\left(\frac{v_t}{v}-\frac{|\nabla
v|^2}{v}\right)+\alpha\frac{v_t}{v}\frac{|\nabla
v|^2}{v}\\
=&(2\alpha+2)\frac{v_t}{v}\frac{|\nabla v|^2}{v}-3\frac{|\nabla
v|^4}{v^2}-\alpha\left(\frac{v_t}{v}\right)^2.
\end{alignat}
Therefore, \eqref{sec6} and \eqref{sec7} give
\begin{alignat}{1}\label{sec8}
2(m-1)&v\nabla\left(\frac{|\nabla v|^2}{v}\right)\nabla(\log
v)-2\alpha(m-1)v\nabla\left(\frac{v_t}{v}\right)\nabla(\log
v)+\frac{2}{v}\nabla v\nabla|\nabla v|^2\\
\zdhy&-\alpha\frac{2}{v}\nabla v\nabla v_t+2(m-1)\Delta
v\frac{|\nabla v|^2}{v}-\frac{|\nabla v|^4}{v^2} -\alpha(m-1)\Delta
v\frac{v_t}{v}+\alpha\frac{v_t}{v}\frac{|\nabla v|^2}{v}\\
\zdhy=&2m \nabla v\nabla F-\left(\frac{v_t}{v}- \frac{|\nabla
v|^2}{v}\right)^2+(1-\alpha)\left(\frac{v_t}{v}\right)^2\\ \zdhy
=&2m \nabla v\nabla F-((m-1)\Delta
v)^2+(1-\alpha)\left(\frac{v_t}{v}\right)^2.
\end{alignat}
Putting \eqref{sec8} into \eqref{sec5} yields $$\aligned
\mathcal{L}(F)=&-2(m-1)v_{ij}^2-2(m-1)R_{ij}v_iv_j+2m \nabla v\nabla F\\
&-((m-1)\Delta
v)^2+(1-\alpha)(\frac{v_t}{v})^2-\alpha'\frac{v_t}{v}-\varphi'.
\endaligned$$ It completes the proof of Lemma 2.1.

Now we prove Theorem 1.1. Define $\tilde{F}=\frac{|\nabla v|^2}{v}-
\alpha \frac{v_t}{v}$, where $\alpha>1$ is a constant. Then we have
from \eqref{sec9}
$$\mathcal{L}(\tilde{F})=-2(m-1)v_{ij}^2-2(m-1)R_{ij}v_iv_j
+2m \nabla v \nabla \tilde{F}-((m-1)\Delta v)^2 +(1-\alpha)\left(
\frac{v_t}{v} \right)^2.$$ Under the assumption that $\mathrm{Ric}
\geq - K$ and the definition of $M$, we have
\begin{equation}\label{2sechuangsec1}\aligned
\mathcal{L}(\tilde{F})=&-2(m-1)v_{ij}^2-2(m-1)R_{ij}v_iv_j +2m
\nabla v \nabla \tilde{F}
-\left( (m-1)\Delta v \right)^2 +(1-\alpha)\left( \frac{v_t}{v} \right)^2\\
\leq& -\frac{2}{n (m-1)}((m-1)\Delta v)^2+2(m-1) K |\nabla v|^2 +2m
\nabla v \nabla \tilde{F}
-\left( (m-1)\Delta v \right)^2\\
\leq& -\frac{1}{a}((m-1) \Delta v)^2+2(m-1)MK\frac{|\nabla v|^2}{v}
+2 m \nabla v \nabla \tilde{F}\\
=& -\frac{1}{a\alpha^2}\left(\tilde{F}+(\alpha-1)\frac{|\nabla
v|^2}{v}\right)^2+2(m-1)MK\frac{|\nabla v|^2}{v} +2 m \nabla v
\nabla \tilde{F},
\endaligned\end{equation}
 where the last equality used
$$(m-1)\Delta v=\frac{v_t}{v}-\frac{|\nabla
v|^2}{v}=-\frac{1}{\alpha}\left(\tilde{F}+(\alpha-1)\frac{|\nabla
v|^2}{v}\right).$$

Denote by $B_p(R)$ the geodesic ball centered at $p$ with radius
$R$. Take a cut-off function $\phi$ (see \cite{schoenyau})
satisfying $\mathrm{supp}(\phi) \subset B_p(2R)$, $\phi |_{B_p(R)}
=1$ and
\begin{gather}
\frac{|\nabla \phi|^2}{\phi}\leq \frac{C}{R^2}, \\
-\Delta \phi \leq \frac{C}{R^2}\left( 1+ \sqrt{K} \coth (\sqrt{K}
R)\right),\label{phi-prop}
\end{gather}
where $C$ is a constant depending only on $n$. Define
$G=t\phi \tilde{F}$. Next we will apply maximum principle to $G$ on
$B_p(2 R)\times [0,T]$. Assume $G$ achieves its maximum at the point
$(x_0,s)\in B_p(2 R)\times [0,T]$ and assume $G(x_0,s)>0$
(otherwise the proof is trivial), which implies $s>0$. Then at the
point $(x_0,s)$, it holds that
$$\mathcal{L}(G)\geq 0, \ \ \
\nabla \tilde{F}=-\frac{\tilde{F}}{\phi} \nabla \phi$$
and by use of \eqref{2sechuangsec1}, we have
\begin{alignat*}{1}
0\leq&
\mathcal{L}(G)=s\phi\mathcal{L}(\tilde{F})
-s(m-1)v\tilde{F}\Delta\phi-2s(m-1)v\nabla
\tilde{F} \nabla\phi+\phi \tilde{F}\\ \zdhy
=&s\phi\mathcal{L}(\tilde{F})-(m-1)v\frac{\Delta\phi}{\phi}G+2(m-1)
v\frac{|\nabla\phi|^2}{\phi^2}G
 +\frac{G}{s}\\ \zdhy
\leq&s\phi\left(-\frac{1}{a\alpha^2}(\tilde{F}+(\alpha-1)\frac{|\nabla
v|^2}{v})^2+2(m-1)MK\frac{|\nabla v|^2}{v} -2m\nabla v\frac{\nabla
\phi}{\phi}\tilde{F} \right) \\ \zdhy &-(m-1)v\frac{\Delta
\phi}{\phi}G+2(m-1)v\frac{|\nabla \phi|^2}{\phi^2}G+\frac{G}{s}.
\end{alignat*}
Let $\frac{|\nabla v |^2}{v}=\mu \tilde{F}$ at the point $(x_0,s)$. Then we have $\mu\geq0$
and \begin{alignat}{1}\label{2sechuangsec2} 0\leq&-\frac{1}{a
\alpha^2} s \phi \tilde{F}^2(1+(\alpha-1)\mu)^2+2(m-1) MK s \phi \mu
\tilde{F}
-2ms\phi \nabla v \frac{\nabla \phi}{\phi}  \tilde{F}  \\
& -(m-1)v   \frac{\Delta \phi}{\phi}G
+2(m-1)v\frac{|\nabla \phi|^2}{\phi^2}G+\frac{G}{s} \\
\leq&-\frac{1}{a \alpha^2 s \phi } (1+(\alpha-1)\mu)^2G^2+2(m-1) MK
\mu G +2m G\frac{|\nabla
\phi|}{s^{\frac{1}{2}}\phi^{\frac{3}{2}}}M^{\frac{1}{2}}
\mu^{\frac{1}{2}} G^{\frac{1}{2}}\\
&-(m-1)v   \frac{\Delta \phi}{\phi}G +2(m-1)v\frac{|\nabla
\phi|^2}{\phi^2}G+\frac{G}{s}.
\end{alignat}
Multiplying the both sides of \eqref{2sechuangsec2} with
$\frac{\phi}{G}$ yields
\begin{alignat}{1}\label{2sechuangsec3}
\frac{1}{a \alpha^2 s}& (1+(\alpha-1)\mu)^2 G-  2  m \frac{|\nabla
\phi|}{s^{\frac{1}{2}}\phi^{\frac{1}{2}}}
M^{\frac{1}{2}} \mu^{\frac{1}{2}} G^{\frac{1}{2}} \\
 &\quad \leq  2(m-1) MK \mu \phi-(m-1)v  \Delta
 \phi + 2(m-1)v \frac{|\nabla \phi|^2}{\phi} +\frac{\phi}{s} \\
&\quad \leq  2(m-1) MK \mu -(m-1)M\Delta \phi + 2(m-1)M
\frac{|\nabla \phi|^2}{\phi} +\frac{1}{s}.
\end{alignat} For the inequality $Ax^2-2Bx\leq C$, we have
$x\leq\frac{2B}{A}+\left(\frac{C}{A}\right)^{\frac{1}{2}}$. Applying
this inequality to \eqref{2sechuangsec3} by setting
$x=G^{\frac{1}{2}}$ gives
\begin{alignat}{1}\label{2sechuangsec4}
G^{\frac{1}{2}}\leq&\frac{2a\alpha^2mM^{\frac{1}{2}} s^{\frac{1}{2}}
 \mu^{\frac{1}{2}}}{(1+(\alpha-1)\mu)^2}
\frac{|\nabla \phi|}{\phi^{\frac{1}{2}}}+\left\{ \frac{a \alpha^2
s}{(1+(\alpha-1)\mu)^2}
\left( \frac{1}{s}+2(m-1) MK \mu\right.\right.\\
&\left.\left.-(m-1)M\Delta \phi+2(m-1)M\frac{|\nabla \phi|^2}{\phi}
\right)\right\}^{\frac{1}{2}}.
\end{alignat}

Obviously, $$\aligned \frac{2   a \alpha^2m
M^{\frac{1}{2}}s^{\frac{1}{2}}
\mu^{\frac{1}{2}}}{(1+(\alpha-1)\mu)^2} \frac{|\nabla
\phi|}{\phi^{\frac{1}{2}}}=& \frac{2   a \alpha^2mM^{\frac{1}{2}}
s^{\frac{1}{2}}
((\alpha-1)\mu)^{\frac{1}{2}}}{(\alpha-1)^{\frac{1}{2}}(1+(\alpha-1)\mu)^2}
\frac{|\nabla \phi|}{\phi^{\frac{1}{2}}} \\ \zdhy \leq&\frac{ a
\alpha^2 mM^{\frac{1}{2}}s^{\frac{1}{2}}
}{(\alpha-1)^{\frac{1}{2}}}\frac{|\nabla \phi|}{\phi^{\frac{1}{2}}}
\endaligned$$
and \begin{alignat*}{1}\frac{a \alpha^2
s}{(1+(\alpha-1)\mu)^2}\left( \frac{1}{s}+2(m-1)MK\mu\right)\leq
&a\alpha^2+\frac{2a\alpha^2s (m-1)
M K}{\alpha-1}\frac{(\alpha-1)\mu}{(1+(\alpha-1)\mu)^2}\\
\leq&a\alpha^2\left(1+\frac{(m-1)MKs}{2(\alpha-1)}\right).
\end{alignat*}
Thus, by use of \eqref{phi-prop}, we have $$\aligned G^{\frac{1}{2}}(x,T)\leq
G^{\frac{1}{2}}(x_0,s)\leq & \frac{ a \alpha^2 mM^{\frac{1}{2}}s^{\frac{1}{2}}
 } {(\alpha-1)^{\frac{1}{2}}}\frac{|\nabla
\phi|}{\phi^{\frac{1}{2}}} + \left\{a \alpha^2\left( 1 +\frac{(m-1)
M K s}{2(\alpha -1)}
\right.\right.\\
& \left. \left. -(m-1)Ms\Delta \phi + 2(m-1)Ms
\frac{|\nabla \phi|^2}{\phi} \right)\right\}^{\frac{1}{2}} \\
\leq&  \frac{ a \alpha^2 m M^{\frac{1}{2}}T^{\frac{1}{2}}
}{(\alpha-1)^{\frac{1}{2}}}\frac{C}{R}+a^{\frac{1}{2}}T^{\frac{1}{2}}\alpha
\left\{
\frac{1}{T} +\frac{(m-1)MK}{2(\alpha -1)} \right. \\
&\left. +(m-1)M\frac{C}{R^2}\left(1+ \sqrt{K} \coth (\sqrt{K}
R)\right) \right\}^{\frac{1}{2}}.
\endaligned$$
Hence, for all $x \in B_p(R)$, it holds that
\begin{equation}\label{2sechuangsec6}\aligned
\tilde{F}^{\frac{1}{2}}(x,T)\leq&\frac{ a \alpha^2m M^{\frac{1}{2}}
}{(\alpha-1)^{\frac{1}{2}}}\frac{C}{R}+a^{\frac{1}{2}}\alpha \left\{
\frac{1}{T} +\frac{(m-1)MK}{2(\alpha -1)} \right. \\
&\left. +(m-1)M\frac{C}{R^2}\left(1+ \sqrt{K} \coth (\sqrt{K}
R)\right) \right\}^{\frac{1}{2}}.
\endaligned\end{equation} Since $T$ is arbitrary, we obtain, for $x\in B_p(R)$
$$\aligned \tilde{F}^{\frac{1}{2}}(x,t)\leq&\frac{ a \alpha^2m
M^{\frac{1}{2}}
}{(\alpha-1)^{\frac{1}{2}}}\frac{C}{R}+a^{\frac{1}{2}}\alpha \left\{
\frac{1}{t} +\frac{(m-1)MK}{2(\alpha -1)}\right. \\
&\left.+(m-1)M\frac{C}{R^2}\left(1+\sqrt{K}\coth(\sqrt{K} R)\right)
\right\}^{\frac{1}{2}}.
\endaligned$$
We complete the proof of Theorem 1.1.

\vspace*{1mm}

\noindent{\bf Proof of Corollary 1.2.} Along the line of Li-Yau, we
will establish Harnack inequality from a general estimate
\begin{equation}\label{general-est}
\frac{|\nabla v|^2}{v}-\alpha(t)\frac{v_t}{v}-\varphi(t) \leq 0.
\end{equation}
Rewrite \eqref{general-est} as
$$-\frac{v_t}{v}\leq\frac{1}{\alpha(t)}\left(\varphi(t)
-\frac{|\nabla v|^2}{v}\right).$$ Let $f=\log v$. Then we have
$$\aligned
-f_t=&-\frac{v_t}{v}\leq\frac{1}{\alpha(t)}\left(\varphi(t)-\frac{|\nabla
v|^2}{v}\right)\\
\leq&\frac{1}{\alpha(t)}(\varphi(t)-\tilde{M}|\nabla f|^2
).\endaligned$$ Let $\gamma $ be a shortest geodesic joining $x_1$
and $x_2$, $\gamma:[t_1,t_2]\rightarrow M^n$, $\gamma(t_1)=x_1$,
$\gamma(t_2)=x_2$. We define a curve $\zeta$ in $M^n\times (0,\infty)$,
$\zeta :[t_1,t_2]\rightarrow M^n\times (0,\infty)$ by $\zeta(t)=\big(
\gamma(t), t \big)$. Then we have $\zeta(t_1)=(x_1,t_1)$,
$\zeta(t_2)=(x_2,t_2)$. Denote by $\rho=d(x_1,x_2)$, then we have
$|\dot \gamma|=\frac{\rho}{(t_2-t_1)}$ and
\begin{alignat}{1}\label{general-Har}
f(x_1,t_1)-f(x_2,t_2)=&\int_{t_2}^{t_1}\frac{d}{dt}f(\zeta(t))dt\\
\zdhy =&\int_{t_2}^{t_1}\big(\langle\dot \gamma,\nabla
f\rangle+f_t\big)dt\\
\zdhy=&\int_{t_1}^{t_2}\big(-\langle\dot \gamma,\nabla f\rangle+(-f_t)\big)dt\\
\zdhy\leq& \int_{t_1}^{t_2}\left(|\dot \gamma||\nabla
f|+\frac{1}{\alpha(t)}(\varphi(t)-\tilde{M}|\nabla f|^2)\right) dt \\
=&\int_{t_1}^{t_2}\left(-\frac{\tilde{M}}{\alpha(t)}|\nabla f|^2 +|\dot \gamma| |\nabla f| \right)dt +\int_{t_1}^{t_2}\frac{\varphi(t)}{\alpha(t)}dt\\
\leq &\frac{\rho^2}{4\tilde{M}(t_2-t_1)^2}\int_{t_1}^{t_2} \alpha(t)
dt +\int_{t_1}^{t_2}\frac{\varphi(t)}{\alpha(t)}dt,
\end{alignat}
where the last inequality used $-Ax^2+Bx \leq \frac{B^2}{4A}$ and
$|\dot \gamma|=\frac{\rho}{(t_2-t_1)}$.

Let $\alpha>1$ be a constant,
$\varphi=\frac{\alpha^2}{2(\alpha-1)}a(m-1)MK+ \frac{a\alpha^2}{t}$.
We have from \eqref{general-Har}
\begin{equation}\label{sec17}\aligned
f(x_1,t_1)-f(x_2,t_2) \leq&\int_{t_1}^{t_2}
\left\{\frac{\alpha\rho^2}{4\tilde{M}(t_2-t_1)^2}
+\left(\frac{\alpha}{2(\alpha-1)}a(m-1)MK+
\frac{a\alpha}{t}\right)\right\} dt\\
=&\frac{\alpha\rho^2}{4\tilde{M}(t_2-t_1)}
+\frac{\alpha}{2(\alpha-1)}a(m-1)MK(t_2-t_1)+a\alpha\log\frac{t_2}{t_1}.
\endaligned\end{equation}
Therefore, we arrive at
$$v(x_1,t_1)\leq v(x_2,t_2)\left(\frac{t_2}{t_1}
\right)^{a\alpha}\exp
\left\{\frac{\alpha\rho^2}{4\tilde{M}(t_2-t_1)}
+\frac{\alpha}{2(\alpha-1)}a(m-1)MK(t_2-t_1) \right\}.$$ We complete the proof of Corollary 1.2.

\section{Proof of Theorem 1.2}

Define $\overline{F}=\frac{|\nabla v|^2}{v}-\alpha\frac{v_t}{v}$,
where $\alpha=e^{2(m-1)MKt}$ is a function depending on $t$. Under
the assumption that ${\rm Ric}\geq-K$, we have from \eqref{sec9}
\begin{alignat}{1}\label{hamilton1}
\mathcal{L}(\overline{F})\leq&-\frac{2}{n(m-1)}((m-1)\Delta
v)^2+2(m-1)MK \frac{|\nabla v|^2}{v}+2m \nabla v\nabla
\overline{F}\\ \zdhy &-((m-1)\Delta
v)^2+(1-\alpha)\left(\frac{v_t}{v}\right)^2-\alpha'\frac{v_t}{v}\\
\zdhy =&-\frac{1}{a}((m-1)\Delta v)^2+2(m-1)MK\frac{|\nabla
v|^2}{v}+2m \nabla v\nabla \overline{F}-\alpha'\frac{v_t}{v}\\ \zdhy
&+(1-\alpha)\left(\frac{v_t}{v}\right)^2,
\end{alignat} and hence
\begin{alignat}{2}\label{hamilton2}
\mathcal{L}(\alpha^{-1}\overline{F})=&(\alpha^{-1})'\overline{F}
+\alpha^{-1}\mathcal{L}(\overline{F})\\ \zdhy
\leq&(\alpha^{-1})'\frac{|\nabla
v|^2}{v}-\alpha(\alpha^{-1})'\frac{v_t}{v}\\
\zdhy &-\frac{1}{a}\alpha^{-1}((m-1)\Delta
v)^2+2(m-1)MK\alpha^{-1}\frac{|\nabla v|^2}{v}+2m \alpha^{-1}\nabla
v\nabla \overline{F}-\alpha'\alpha^{-1}\frac{v_t}{v}\\ \zdhy
&+(1-\alpha)\alpha^{-1}\left(\frac{v_t}{v}\right)^2\\ \zdhy
=&-\frac{1}{a}\alpha^{-1}((m-1)\Delta v)^2+2m \alpha^{-1}\nabla
v\nabla
\overline{F}+(1-\alpha)\alpha^{-1}\left(\frac{v_t}{v}\right)^2\\
\zdhy \leq&-\frac{1}{a}\alpha^{-1}((m-1)\Delta v)^2+2m
\alpha^{-1}\nabla v\nabla \overline{F}\\ \zdhy
=&-\frac{1}{a}\alpha^{-1}\left((\alpha^{-1}-1)\frac{|\nabla
v|^2}{v}-\alpha^{-1}\overline{F}\right)^2+2m \nabla v\nabla
(\alpha^{-1}\overline{F}),
\end{alignat}
where we used $(m-1)\Delta v=(\alpha^{-1}-1)\frac{|\nabla
v|^2}{v}-\alpha^{-1}\overline{F}$.

Recall that one can construct a cut-off function $\phi$ as before,
which satisfies ${\rm supp}(\phi)\subset B_p(2R)$,
$\phi|_{B_p(R)}=1$ and
$$\frac{|\nabla\phi|^2}{\phi}\leq\frac{C}{R^2},$$
$$-\Delta\phi\leq\frac{C}{R^2}\left(1+\sqrt{K}R\coth(\sqrt{K}R)\right),$$
where $C$ is a constant depending only on $n$. Define
$G=t\phi\alpha^{-1}\overline{F}$. Next we are to apply the maximum
principle to $G$ on $B_p(2R)\times [0,T]$. Assume $G$ achieves its
maximum at the point $(x_0,s)\in B_p(2R)\times [0,T]$, and assume
that $G(x_0,s)>0$ (otherwise the proof is trivial), which implies
$s>0$. Then at the point $(x_0,s)$, it holds that
$$\mathcal{L}(G)\geq0,\ \ \ \ \nabla
(\alpha^{-1}\overline{F})=-\frac{\alpha^{-1}\overline{F}}{\phi}\nabla\phi$$
and
\begin{alignat}{1}\label{hamilton3}
0\leq& \mathcal{L}(G)=s\phi\mathcal{L}(\alpha^{-1}\overline{F})
-s(m-1)v\alpha^{-1}\overline{F}\Delta\phi-2s(m-1)v\nabla
(\alpha^{-1}\overline{F}) \nabla\phi+\phi \alpha^{-1}\overline{F}\\
\zdhy
=&s\phi\mathcal{L}(\alpha^{-1}\overline{F})-(m-1)v\frac{\Delta\phi}{\phi}G+2(m-1)
v\frac{|\nabla\phi|^2}{\phi^2}G
 +\frac{G}{s}\\ \zdhy
\leq&-\frac{s\phi\alpha^{-1}}{a}\left((\alpha^{-1}-1)\frac{|\nabla
v|^2}{v}-\alpha^{-1}\overline{F}\right)^2+2m s\phi\nabla
v\nabla (\alpha^{-1}\overline{F}) \\
\zdhy &-(m-1)v\frac{\Delta \phi}{\phi}G+2(m-1)v\frac{|\nabla
\phi|^2}{\phi^2}G+\frac{G}{s}\\ \zdhy
\leq&-\frac{s\phi\alpha^{-1}}{a}(\alpha^{-1}-1)^2\frac{|\nabla
v|^4}{v^2}-\frac{\alpha^{-1}}{as\phi}G^2+\frac{2\alpha^{-1}(\alpha^{-1}-1)}{a}\frac{|\nabla
v|^2}{v}G-2m \nabla
v\frac{\nabla \phi}{\phi}G \\
\zdhy &-(m-1)v\frac{\Delta \phi}{\phi}G+2(m-1)v\frac{|\nabla
\phi|^2}{\phi^2}G+\frac{G}{s}\\ \zdhy
\leq&-\frac{\alpha^{-1}}{as\phi}G^2+\frac{2\alpha^{-1}(\alpha^{-1}-1)}{a}\frac{|\nabla
v|^2}{v}G+2m M^{\frac{1}{2}}\frac{|\nabla \phi|}{\phi}\frac{|\nabla
v|}{v^{\frac{1}{2}}}G \\ \zdhy &-(m-1)M\frac{\Delta
\phi}{\phi}G+2(m-1)M\frac{|\nabla \phi|^2}{\phi^2}G+\frac{G}{s}.
\end{alignat}
Multiplying the both sides of \eqref{hamilton3} with
$\frac{\alpha as\phi}{G}$ yields
\begin{alignat}{1}\label{hamilton4}
G(x,T)\leq& G(x_0,s)\leq-2(1-\alpha^{-1})s\phi\frac{|\nabla
v|^2}{v}+2mM^{\frac{1}{2}}as\alpha|\nabla \phi|\frac{|\nabla
v|}{v^{\frac{1}{2}}}-(m-1)Mas\alpha\Delta \phi\\
\zdhy &+2(m-1)Mas\alpha\frac{|\nabla \phi|^2}{\phi}+a\alpha\phi\\
\zdhy \leq&\frac{m^2Ma^2\alpha^2s}{2(1-\alpha^{-1})}\frac{|\nabla
\phi|^2}{\phi}-(m-1)Mas\alpha\Delta
\phi+2(m-1)Mas\alpha\frac{|\nabla \phi|^2}{\phi}+a\alpha\phi\\
\zdhy
\leq&\left(\frac{m^2Ma^2\alpha^2T}{2(1-\alpha^{-1})}+3(m-1)MaT\alpha\right)\frac{C}{R^2}\\&
+(m-1)MaT\alpha\sqrt{K}\coth(\sqrt{K}R)\frac{C}{R}+a\alpha.
\end{alignat}
Hence, for all $x \in B_p(R)$, it holds that
\begin{alignat*}{1}(\alpha^{-1}\overline{F})(x,T)\leq&
\left(\frac{m^2Ma^2\alpha^2}{2(1-\alpha^{-1})}+3(m-1)Ma\alpha\right)\frac{C}{R^2}\\&
+(m-1)Ma\alpha\sqrt{K}\coth(\sqrt{K}R)\frac{C}{R}+\frac{a\alpha}{T},
\end{alignat*}
and hence,
\begin{alignat*}{1}\overline{F}(x,T)\leq&
\left(\frac{m^2Ma^2\alpha^3}{2(1-\alpha^{-1})}+3(m-1)Ma\alpha^2\right)\frac{C}{R^2}\\&
+(m-1)Ma\alpha^2\sqrt{K}\coth(\sqrt{K}R)\frac{C}{R}+\frac{a\alpha^2}{T},
\end{alignat*}
Since $T$ is arbitrary, we complete the proof of Theorem 1.2.

\vspace*{1mm}

\noindent{\bf Proof of Corollary 1.4.} Choosing  $\alpha(t)=e^{2(m-1)MKt},\varphi(t)=\frac{a\alpha^2(t)}{t}$ in \eqref{general-Har}, we get
\begin{equation}\label{sec17}\aligned
\log v(x_1,t_1)-\log v(x_2,t_2)
\leq&\int_{t_1}^{t_2}\left( \frac{\rho^2\alpha}{4\tilde{M}(t_2-t_1)^2}
 +\frac{a\alpha}{t}\right)dt\\
\leq&\int_{t_1}^{t_2}\left(
\frac{\rho^2\alpha}{4\tilde{M}(t_2-t_1)^2}
+\frac{a\alpha}{t_1}\right)dt\\
=&\left( \frac{\rho^2}{4\tilde{M}(t_2-t_1)^2}
+\frac{a}{t_1}\right)\frac{e^{2(m-1)MKt_2}-e^{2(m-1)MKt_1}}{2(m-1)MK},
\endaligned\end{equation} which concludes the proof of Corollary
1.4.

\section{Proof of Theorem 1.3}

Under the assumption that ${\rm Ric}\geq-K$, we have from
\eqref{sec9}
\begin{alignat}{1}\label{sec10}
\mathcal{L}(F)\leq&-\frac{1}{a}((m-1)\Delta v)^2+2(m-1)MK\frac{|\nabla
v|^2}{v}+2m \nabla v\nabla F-\alpha'\frac{v_t}{v}\\ \zdhy
&-\varphi'+(1-\alpha)\left(\frac{v_t}{v}\right)^2\\ \zdhy
=&-\frac{1}{a}((m-1)\Delta
v+\varphi)^2+\frac{2}{a}\varphi((m-1)\Delta v
)+\frac{1}{a}\varphi^2+2(m-1)MK\frac{|\nabla v|^2}{v}\\ \zdhy
&+2m\nabla v\nabla
F-\alpha'\frac{v_t}{v}-\varphi'+(1-\alpha)\left(\frac{v_t}{v}\right)^2\\
\zdhy =&-\frac{1}{a}((m-1)\Delta
v+\varphi)^2+\frac{2}{a}\varphi\left(\frac{v_t}{v}-\frac{|\nabla
v|^2}{v}\right)+\frac{1}{a}\varphi^2+2(m-1)MK\frac{|\nabla v|^2}{v}\\
\zdhy &+2m\nabla v\nabla
F-\alpha'\frac{v_t}{v}-\varphi'+(1-\alpha)\left(\frac{v_t}{v}\right)^2\\
\zdhy =&-\frac{1}{a}((m-1)\Delta v+\varphi)^2+2m\nabla v\nabla F\\
\zdhy &-\left(\frac{2}{a}\varphi-2(m-1)MK\right)\left(\frac{|\nabla
v|^2}{v}-\frac{\frac{2}{a}\varphi-\alpha'}{\frac{2}{a}\varphi-2(m-1)MK}
\frac{v_t}{v}-\varphi\right)\\ \zdhy
&-\left(\frac{2}{a}\varphi-2(m-1)MK\right)\varphi+\frac{1}{a}\varphi^2-\varphi'
+(1-\alpha)\left(\frac{v_t}{v}\right)^2.
\end{alignat}
Take \begin{equation}\label{addsec11}\aligned \varphi(t)
=&a(m-1)MK\{\coth((m-1)MKt)+1\}\\
\alpha(t)=&1+\frac{\cosh((m-1)MKt)\sinh((m-1)MKt)-(m-1)MKt}{\sinh^2((m-1)MKt)},
\endaligned\end{equation} then $\alpha(t)$ and $\varphi(t)$ satisfy
 the following equations:
\begin{equation}\label{addsec12}\left\{\begin{array}{l}
-\left(\frac{2}{a}\varphi-2(m-1)MK\right)\varphi+\frac{1}{a}\varphi^2-\varphi'=0\\
\frac{\frac{2}{a}\varphi-\alpha'}{\frac{2}{a}\varphi-2(m-1)MK}=\alpha.
\end{array}\right.\end{equation}
Moreover, it is easy to see that $\alpha\geq1$. Putting
\eqref{addsec12} into \eqref{sec10}, we obtain
\begin{equation}\label{sec11}\aligned\mathcal{L}(F) \leq&
-\frac{1}{a}((m-1)\Delta v+\varphi)^2+2m\nabla v\nabla
F-2(m-1)MK\coth((m-1)MKt)F\\
=&-\frac{1}{a\alpha^2}\left(F+(\alpha-1)\left(\frac{|\nabla
v|^2}{v}-\varphi\right)\right)^2+2m\nabla v\nabla
F\\
&-2(m-1)MK\coth((m-1)MKt)F,
\endaligned\end{equation} where we used
$$(m-1)\Delta v+\varphi=\frac{v_t}{v}-\frac{|\nabla
v|^2}{v}+\varphi=-\frac{1}{\alpha}\left(F+(\alpha-1)\left(\frac{|\nabla
v|^2}{v}-\varphi\right)\right).$$

We can construct a cut-off function $\phi$ as before,
satisfying  ${\rm supp}(\phi)\subset B_p(2R)$,
$\phi|_{B_p(R)}=1$ and
$$\frac{|\nabla\phi|^2}{\phi}\leq\frac{C}{R^2},$$
$$-\Delta\phi\leq\frac{C}{R^2}\left(1+\sqrt{K}R\coth(\sqrt{K}R)\right),$$
where $C$ is a constant which depends only on $n$. Define $G=\beta(t)\phi F$, where
$\beta(t)$ is a positive function to be determined. Next we are to
apply the maximum principle to $G$ on $B_p(2R)\times [0,T]$. Assume
$G$ achieves its maximum at the point $(x_0,s)\in B_p(2R)\times
[0,T]$, and assume that $G(x_0,s)>0$ (otherwise the proof is
trivial), which implies $s>0$. Then at the point $(x_0,s)$, it
holds that
$$\mathcal{L}(G)\geq0,\ \ \ \ \nabla
F=-\frac{F}{\phi}\nabla\phi$$ and
\begin{alignat}{1}\label{sec12}
0\leq&\mathcal{L}(G)=\beta'\phi F+\beta\phi
\mathcal{L}(F)-(m-1)\beta vF\Delta \phi-2(m-1)\beta
v\nabla\phi\nabla F\\ \zdhy =&\frac{\beta'}{\beta}G+\beta\phi
\mathcal{L}(F)-(m-1)v\frac{\Delta
\phi}{\phi}G+2(m-1)v\frac{|\nabla\phi|^2}{\phi^2}G\\ \zdhy
\leq&\frac{\beta'}{\beta}G-\frac{\beta\phi}{a\alpha^2}
\left(F+(\alpha-1)(\frac{|\nabla
v|^2}{v}-\varphi)\right)^2+2m\beta\phi\nabla v\nabla F\\ \zdhy
&-2(m-1)\beta\phi MK\coth((m-1)MKs)F-(m-1)v\frac{\Delta
\phi}{\phi}G+2(m-1)v\frac{|\nabla\phi|^2}{\phi^2}G\\ \zdhy
=&\frac{\beta'}{\beta}G-\frac{1}{a\alpha^2\beta\phi}G^2
-\frac{\beta\phi(\alpha-1)^2}{a\alpha^2}\left(\frac{|\nabla
v|^2}{v}-\varphi\right)^2-2\frac{(\alpha-1)}{a\alpha^2}\left(\frac{|\nabla
v|^2}{v}-\varphi\right)G-2m\nabla v\frac{\nabla \phi}{\phi}G\\ \zdhy
&-2(m-1)MK\coth((m-1)MKs)G-(m-1)v\frac{\Delta
\phi}{\phi}G+2(m-1)v\frac{|\nabla\phi|^2}{\phi^2}G\\ \zdhy
\leq&\frac{\beta'}{\beta}G-\frac{1}{a\alpha^2\beta\phi}G^2
-2\frac{(\alpha-1)}{a\alpha^2}\left(\frac{|\nabla
v|^2}{v}-\varphi\right)G-2m\nabla v\frac{\nabla \phi}{\phi}G\\ \zdhy
&-2(m-1)MK\coth((m-1)MKs)G-(m-1)v\frac{\Delta
\phi}{\phi}G+2(m-1)v\frac{|\nabla\phi|^2}{\phi^2}G.
\end{alignat} Multiplying the both sides of
\eqref{sec12} with $\frac{a\alpha^2\beta\phi}{G}$ yields
\begin{alignat}{1}\label{sec13}
G(x_0,s)\leq&a\alpha^2\beta\left\{\frac{\beta'}{\beta}\phi
-2\frac{(\alpha-1)}{a\alpha^2}(\frac{|\nabla
v|^2}{v}-\varphi)\phi-2m\nabla v\nabla \phi\right.\\ \zdhy
&\left.-2(m-1)MK\coth((m-1)MKs)\phi-(m-1)v\Delta \phi\right.\\ \zdhy
&\left.+2(m-1)v\frac{|\nabla\phi|^2}{\phi}\right\}\\ \zdhy
=&a\alpha^2\beta\left\{\left(\frac{\beta'}{\beta}
+2\frac{(\alpha-1)}{a\alpha^2}\varphi
-2(m-1)MK\coth((m-1)MKs)\right)\phi\right.\\ \zdhy
&\left.-(m-1)v\Delta
\phi+2(m-1)v\frac{|\nabla\phi|^2}{\phi}\right.\\ \zdhy
&\left.-2\frac{(\alpha-1)\phi}{a\alpha^2}\frac{|\nabla
v|^2}{v}+2mv^{\frac{1}{2}}|\nabla \phi|\frac{|\nabla
v|}{v^{\frac{1}{2}}}\right\}\\ \zdhy
\leq&\beta\left\{2(\alpha-1)\varphi-a\alpha^2\left(
2(m-1)MK\coth((m-1)MKs)-\frac{\beta'}{\beta}\right)\right\}\phi\\
\zdhy &+a\alpha^2\beta \left\{-(m-1)\Delta
\phi+2(m-1)\frac{|\nabla\phi|^2}{\phi}
+\frac{am^2\alpha^2}{2(\alpha-1)}\frac{|\nabla\phi|^2}{\phi}\right\}M,
\end{alignat}
where the last inequality used $-Ax^2+Bx\leq\frac{B^2}{4A}$.
Choosing $\beta(t)=\tanh((m-1)MKt)$, we have
$\frac{\beta'}{\beta}=\frac{(m-1)MK}
{\sinh((m-1)MKt)\cosh((m-1)MKt)}$ and
$$2(\alpha-1)\varphi-a\alpha^2\left(
2(m-1)MK\coth((m-1)MKt)-\frac{\beta'}{\beta}\right)\leq0.$$ Note
that $\alpha,\ \frac{\beta}{\alpha-1}$ is bounded uniformly
and $\beta$ is non-decreasing. Thus from \eqref{sec13} we obtain
\begin{alignat}{1}\label{sec15}
G(x,T)\leq &G (x_0,s)\\ \zdhy \leq&a\alpha^2\beta
\left\{-(m-1)\Delta \phi+2(m-1)\frac{|\nabla\phi|^2}{\phi}
+\frac{am^2\alpha^2}{2(\alpha-1)}\frac{|\nabla\phi|^2}{\phi}\right\}M\\
\zdhy \leq&\left\{a(m-1)\beta(T)
\left(\frac{C}{R^2}+\frac{C\sqrt{K}\coth(\sqrt{K}R)}{R}\right)
+a^2m^2\frac{C}{R^2}\right\}M.
\end{alignat}
Hence, for all $x\in B_p(R)$, one has
$$F(x,T)\leq \left\{a(m-1)
\left(\frac{C}{R^2}+\frac{C\sqrt{K}\coth(\sqrt{K}R)}{R}\right)
+a^2m^2\frac{C}{R^2\beta(T)}\right\}M.$$ Since $T$ is arbitrary, we
obtain
\begin{equation}\label{sec16}\aligned
 \frac{|\nabla v|^2}{v} -
\alpha(t)\frac{v_t}{v}-\varphi(t) \leq&\left\{a(m-1)
\left(\frac{C}{R^2}+\frac{C\sqrt{K}\coth(\sqrt{K}R)}{R}\right)\right.\\
&\left.+a^2m^2\frac{C}{R^2\tanh((m-1)MKt)}\right\}M.
\endaligned\end{equation}
We complete the proof of Theorem 1.3.

\vspace*{1mm}

\noindent{\bf Proof of Corollary 1.6.}
Putting $\alpha(t),\varphi(t)$ given by \eqref{maysec7} into
\eqref{general-Har} gives
\begin{alignat*}{1}
&\log v(x_1,t_1)-\log v(x_2,t_2)\\ \zdhy
\leq&\frac{\rho^2}{4\tilde{M}(t_2-t_1)^2}\int_{t_1}^{t_2} \alpha(t)
dt +\int_{t_1}^{t_2}\frac{\varphi(t)}{\alpha(t)}dt\\ \zdhy
=&\left\{\frac{\rho^2}{4\tilde{M}(t_2-t_1)^2}
\left(t+\frac{(m-1)MKt\coth((m-1)MKt)-1}{ (m-1)MK}\right)\right.\\
\zdhy &\left.+\frac{a}{2}\log\frac{\sinh(2(m-1)MKt)+\cosh(2(m-1)MKt)
-2(m-1)MKt-1}{2(m-1)MK}\right\}\Big|_{t_1}^{t_2}\\
=&\frac{\rho^2}{4\tilde{M}(t_2-t_1)}\left(1+A_2(t_1,t_2)\right)+\log
A_1(t_1,t_2),
\end{alignat*}
where
$$A_1(t_1,t_2)=\left(
\frac{\exp(2(m-1)MKt_2)-2(m-1)MKt_2-1}{\exp(2(m-1)MKt_1)
-2(m-1)MKt_1-1}\right)^{\frac{a}{2}}$$
$$A_2(t_1,t_2)=\frac{t_2\coth((m-1)MKt_2)-t_1\coth((m-1)MKt_1)}{t_2-t_1}.$$
Therefore, we arrive at
$$v(x_1,t_1)\leq v(x_2,t_2)A_1(t_1,t_2)
\exp\left\{\frac{\rho^2}{4\tilde{M}(t_2-t_1)}(1+A_2(t_1,t_2))\right\}.$$
We complete the proof of Corollary 1.6.

\section{Proof of Theorem 1.4}

The proof of Theorem 1.4 is similar to the proof of Theorem 1.3.
Under the assumption that ${\rm Ric}\geq-K$, we have from
\eqref{sec9}
\begin{alignat}{1}\label{proof1}
\mathcal{L}(F)\leq&-\frac{1}{a}((m-1)\Delta v)^2+2(m-1)MK \frac{|\nabla
v|^2}{v}+2m \nabla v\nabla F-\alpha'\frac{v_t}{v}\\ \zdhy
&-\varphi'+(1-\alpha)\left(\frac{v_t}{v}\right)^2\\ \zdhy
=&-\frac{1}{a}\left((m-1)\Delta v+
a\left(\frac{1}{t}+(m-1)MK\right)\right)^2+2(m-1)\Delta
v\left(\frac{1}{t}+(m-1)MK\right)\\ \zdhy
&+a\left(\frac{1}{t}+(m-1)MK\right)^2 +2(m-1)MK\frac{|\nabla
v|^2}{v}+2m\nabla v\nabla F-\alpha'\frac{v_t}{v}\\ \zdhy
&-\varphi'+(1-\alpha)\left(\frac{v_t}{v}\right)^2\\ \zdhy
=&-\frac{1}{a}\left((m-1)\Delta v+
a\left(\frac{1}{t}+(m-1)MK\right)\right)^2-\frac{2}{t}\left\{\frac{|\nabla
v|^2}{v}\right.\\ \zdhy
&\left.-\frac{2(\frac{1}{t}+(m-1)MK)-\alpha'}
{\frac{2}{t}}\frac{v_t}{v}-\varphi\right\}
-\frac{2}{t}\varphi+a\left(\frac{1}{t}+(m-1)MK\right)^2\\ \zdhy
&-\varphi'+(1-\alpha)\left(\frac{v_t}{v}\right)^2+2m\nabla v\nabla F
\end{alignat}
Taking \begin{equation}\label{proof2}\aligned \varphi(t)
=&\frac{a}{t}+a(m-1)MK+\frac{a}{3}((m-1)MK)^2t,\\
\alpha(t)=&1+\frac{2}{3}(m-1)MKt,
\endaligned\end{equation} then $\alpha(t)$ and $\varphi(t)$ satisfy
 the following equations:
\begin{equation}\label{proof3}\left\{\begin{array}{l}
-\frac{2}{t}\varphi+a\left(\frac{1}{t}+(m-1)MK\right)^2
-\varphi'=0\\
\frac{2(\frac{1}{t}+(m-1)MK)-\alpha'}{\frac{2}{t}}=\alpha.
\end{array}\right.\end{equation}
Moreover, it is easy to see that $\alpha\geq1$. Putting
\eqref{proof3} into \eqref{proof1}, we obtain
\begin{alignat}{1}\label{proof4}
 \mathcal{L}(F) \leq&
-\frac{1}{a}\left((m-1)\Delta v+
a\left(\frac{1}{t}+(m-1)MK\right)\right)^2-\frac{2}{t}F+2m\nabla v\nabla F\\
\zdhy =&-\frac{1}{a\alpha^2}\left\{F+(\alpha-1)\frac{|\nabla
v|^2}{v}-\left(\frac{2a}{3}(m-1)MK+\frac{a}{3}((m-1)MK)^2t\right)\right\}^2\\
\zdhy &-\frac{2}{t}F+2m\nabla v\nabla F,
\end{alignat}
where we used
$$\aligned
(m-1)&\Delta v+
a\left(\frac{1}{t}+(m-1)MK\right)\\
&=-\frac{1}{\alpha}\left\{F+(\alpha-1)\frac{|\nabla
v|^2}{v}-\left(\frac{2a}{3}(m-1)MK
+\frac{a}{3}((m-1)MK)^2t\right)\right\}.\endaligned$$

Construct a cut-off function $\phi$ as before,
which satisfies ${\rm supp}(\phi)\subset B_p(2R)$,
$\phi|_{B_p(R)}=1$ and
$$\frac{|\nabla\phi|^2}{\phi}\leq\frac{C}{R^2},$$
$$-\Delta\phi\leq\frac{C}{R^2}\left(1+\sqrt{K}R\coth(\sqrt{K}R)\right),$$
where $C$ is a constant depending only on $n$. Define $G=\beta(t)\phi F$, where
$\beta(t)$ is a positive function to be determined. Assume $G$
achieves its maximum at the point $(x_0,s)\in B_p(2R)\times
[0,T]$, and assume that $G(x_0,s)>0$ (otherwise the proof is
trivial), which implies $s>0$. Then at the point $(x_0,s)$, it
holds that
$$\mathcal{L}(G)\geq0,\ \ \ \ \nabla
F=-\frac{F}{\phi}\nabla\phi$$ and
\begin{alignat}{1}\label{proof5}
0\leq&\mathcal{L}(G)=\beta'\phi F+\beta\phi
\mathcal{L}(F)-(m-1)\beta vF\Delta \phi-2(m-1)\beta
v\nabla\phi\nabla F\\ \zdhy =&\frac{\beta'}{\beta}G+\beta\phi
\mathcal{L}(F)-(m-1)v\frac{\Delta
\phi}{\phi}G+2(m-1)v\frac{|\nabla\phi|^2}{\phi^2}G\\ \zdhy
\leq&\frac{\beta'}{\beta}G-\frac{\beta\phi}{a\alpha^2}
\left\{F+(\alpha-1)\frac{|\nabla
v|^2}{v}-\left(\frac{2a}{3}(m-1)MK+\frac{a}{3}((m-1)MK)^2s\right)\right\}^2\\
\zdhy &-\frac{2}{s}\beta\phi F+2m\beta\phi\nabla v\nabla
F-(m-1)v\frac{\Delta
\phi}{\phi}G+2(m-1)v\frac{|\nabla\phi|^2}{\phi^2}G\\ \zdhy
=&\frac{\beta'}{\beta}G-\frac{G^2}{a\alpha^2\beta\phi}
-\frac{2G}{a\alpha^2}\left\{(\alpha-1)\frac{|\nabla
v|^2}{v}-\left(\frac{2a}{3}(m-1)MK+\frac{a}{3}((m-1)MK)^2s\right)\right\}\\
\zdhy &-\frac{\beta\phi}{a\alpha^2}\left\{(\alpha-1)\frac{|\nabla
v|^2}{v}-\left(\frac{2a}{3}(m-1)MK+\frac{a}{3}((m-1)MK)^2s\right)\right\}^2\\
\zdhy &-\frac{2}{s}G-2m\nabla v\frac{\nabla
\phi}{\phi}G-(m-1)v\frac{\Delta
\phi}{\phi}G+2(m-1)v\frac{|\nabla\phi|^2}{\phi^2}G\\ \zdhy
\leq&\frac{\beta'}{\beta}G-\frac{G^2}{a\alpha^2\beta\phi}
-\frac{2G}{a\alpha^2}\left\{(\alpha-1)\frac{|\nabla
v|^2}{v}-\left(\frac{2a}{3}(m-1)MK+\frac{a}{3}((m-1)MK)^2s\right)\right\}\\
\zdhy &-\frac{2}{s}G-2m\nabla v\frac{\nabla
\phi}{\phi}G-(m-1)v\frac{\Delta
\phi}{\phi}G+2(m-1)v\frac{|\nabla\phi|^2}{\phi^2}G.
\end{alignat}
Multiplying both sides of \eqref{proof5} with
$\frac{a\alpha^2\beta\phi}{G}$ yields
\begin{alignat}{1}\label{proof6}
G(x_0,s)\leq&a\alpha^2\beta\left\{\frac{\beta'}{\beta}\phi
-\frac{2}{a\alpha^2}\left((\alpha-1)\frac{|\nabla
v|^2}{v}-\left(\frac{2a}{3}(m-1)MK+\frac{a}{3}((m-1)MK)^2s\right)\right)\phi\right.\\
\zdhy &\left.-\frac{2}{s}\phi-2m\nabla v\nabla \phi-(m-1)v\Delta
\phi+2(m-1)v\frac{|\nabla\phi|^2}{\phi}\right\}\\ \zdhy
=&a\alpha^2\beta\left\{\left(\frac{\beta'}{\beta}+\frac{2}{3\alpha^2}\Big(2(m-1)MK+
((m-1)MK)^2s\Big)-\frac{2}{s}\right)\phi\right.\\ \zdhy
&\left.-(m-1)v\Delta
\phi+2(m-1)v\frac{|\nabla\phi|^2}{\phi}
-2\frac{(\alpha-1)\phi}{a\alpha^2}\frac{|\nabla
v|^2}{v}+2mv^{\frac{1}{2}}|\nabla \phi|\frac{|\nabla
v|}{v^{\frac{1}{2}}}\right\}\\ \zdhy
\leq&\frac{a\beta}{s}\left\{\frac{2}{3}\Big(2(m-1)MKs+
((m-1)MK)^2s^2\Big)+s\alpha^2\left(\frac{\beta'}{\beta}
-\frac{2}{s}\right)\right\}\phi \\
\zdhy & +a\alpha^2\beta\left\{-(m-1)\Delta
\phi+2(m-1)\frac{|\nabla\phi|^2}{\phi}+\frac{am^2\alpha^2}{2(\alpha-1)}
\frac{|\nabla\phi|^2}{\phi}\right\}M.
\end{alignat}
Choose $\beta(t)=\tanh((m-1)MKt)$, then we have
$\frac{\beta'}{\beta}=\frac{(m-1)MK}
{\sinh((m-1)MKt)\cosh((m-1)MKt)}$ and
$$\frac{2}{3}\Big(2(m-1)MKt+
((m-1)MK)^2t^2\Big)+t\alpha^2\left(\frac{\beta'}{\beta}-\frac{2}{t}\right)\leq0.$$
Note that $\frac{\beta}{\alpha-1}$ is bounded uniformly and
$\alpha,\beta$ is non-decreasing. Thus from \eqref{proof6} we obtain
\begin{alignat*}{1}
G(x,T)\leq &G (x_0,s)\\ \zdhy \leq&a\alpha^2\beta
\left\{-(m-1)\Delta \phi+2(m-1)\frac{|\nabla\phi|^2}{\phi}
+\frac{am^2\alpha^2}{2(\alpha-1)}\frac{|\nabla\phi|^2}{\phi}\right\}M\\
\zdhy \leq&\left\{a(m-1)\alpha^2(T)\beta(T)\left(
\frac{C}{R^2}+\frac{C\sqrt{K}\coth(\sqrt{K}R)}{R}\right)
+a^2m^2\alpha^4(T)\frac{C}{R^2}\right\}M.
\end{alignat*}
Hence, for $x\in B_p(R)$, we have
$$F(x,T)\leq \Big\{a(m-1)\alpha^2(T) \Big(
\frac{C}{R^2}+\frac{C\sqrt{K}\coth(\sqrt{K}R)}{R}\Big)
+\frac{a^2m^2\alpha^4(T)}{\beta(T)}\frac{C}{R^2}\Big\}M.$$ Since $T$
is arbitrary, we obtain
\begin{alignat}{1}\label{proof8}
 \frac{|\nabla v|^2}{v} -
\alpha(t)\frac{v_t}{v}-\varphi(t) \leq&\left\{a(m-1)\alpha^2(t)
\left(\frac{C}{R^2}+\frac{C\sqrt{K}\coth(\sqrt{K}R)}{R}\right)\right.\\
\zdhy
&\left.+\frac{a^2m^2\alpha^4(t)}{\beta(t)}\frac{C}{R^2}\right\}M.
\end{alignat}
It completes the proof of Theorem 1.4.

\vspace*{1mm}

\noindent{\bf Proof of Corollary 1.8.} Recall the estimate \eqref{maysec12}, which implies
$$-\frac{v_t}{v}\leq\frac{1}{\alpha(t)}
\left(\varphi(t)-\frac{|\nabla v|^2}{v}\right),$$ where $\alpha(t)$
and $\varphi(t)$ are given by \eqref{maysec11}. It follows from
\eqref{general-Har} that
\begin{equation}\label{proof9}\aligned
\log v(x_1,t_1)-\log v(x_2,t_2)
\leq\frac{\rho^2}{4\tilde{M}(t_2-t_1)^2}\int_{t_1}^{t_2} \alpha(t)
dt +\int_{t_1}^{t_2}\frac{\varphi(t)}{\alpha(t)}dt.
\endaligned\end{equation}
Choosing $\alpha(t)=1+\frac{2}{3}(m-1)MKt$ and $\varphi(t)
=\frac{a}{t}+a(m-1)MK+\frac{a}{3}((m-1)MK)^2t$ in \eqref{proof9}
concludes the proof of Corollary 1.8.

\noindent {Guangyue Huang} \\
\noindent Department of Mathematical Sciences, Tsinghua University,
Beijing 100084, P.R. China. \ \  E-mail
address:~\textsf{gyhuang$@$math.tsinghua.edu.cn }

\vspace*{1mm}
\noindent {Zhijie Huang} \\
\noindent Department of Mathematical Sciences, Tsinghua University,
Beijing 100084, P.R. China. \ \  E-mail
address:~\textsf{hzj010102$@$126.com }

\vspace*{1mm}
\noindent {Haizhong Li} \\
\noindent Department of Mathematical Sciences, Tsinghua University,
Beijing 100084, P.R. China. \ \  E-mail
address:~\textsf{hli$@$math.tsinghua.edu.cn }

\end{document}